\documentclass[a4paper]{article}

\usepackage{amsfonts,latexsym,amsmath}
\usepackage{graphicx}

\newtheorem{theorem}{Theorem}
\newtheorem{proposition}[theorem]{Proposition}

\newtheorem{corollary}[theorem]{Corollary}
\newtheorem{definition}[theorem]{Definition}

\def\R{\mathbb{R}}

\def\<{\langle}
\def\>{\rangle}
\def\refe#1{(\ref{#1})}
\def\sgn{\mathrm{sgn}}
\def\eps{\varepsilon}
\def\qed{\rule{0.2cm}{0.2cm}}
\def\1{{1\hspace{-1.2mm}{\rm I}}}
\def\limd{\mathop{{\lim}}\limits}
\def\supd{\mathop{{\sup}}\limits}

\begin{document}

\title{A BGK approximation to scalar conservation laws with discontinuous flux}
\author{F. Berthelin\thanks{Laboratoire J. A. Dieudonn\'e, UMR 6621 CNRS,
        Universit\'e de Nice, Parc Valrose,
        06108 Nice cedex 2, France, 
	Email: Florent.Berthelin@unice.fr} and J. Vovelle\thanks{Universit\'e de Lyon ; CNRS ; Universit\'e Lyon 1, Institut Camille Jordan,  43 boulevard du 11 novembre 1918, F-69622 Villeurbanne Cedex, France. Email: vovelle@math.univ-lyon1.fr}}
\date{}
\maketitle

\begin{abstract} We study the BGK approximation to first-order scalar conservation laws with a flux which is discontinuous in the space variable. We show that the Cauchy Problem for the BGK approximation is well-posed and that, as the relaxation parameter tends to $0$, it converges to the (entropy) solution of the limit problem.
\end{abstract}
\medskip

\noindent{\it Keywords}:
scalar conservation laws --  discontinuous flux -- BGK model -- relaxation limit
\medskip

\noindent {\bf Mathematics Subject Classification}:
35L65 -- 35F10 -- 35D05

%\tableofcontents

%%%%%%%%%%%%%%%%%%%%%%%%%%%%%%%%%%%%%%%%%%%%%%%%%%%%%%%%%%%%%%%%%%%%%%%%%%%%%%%%%%%
%%%%%%%%%%%%%%%%%%%%%%%%%%%%%%%%%%%%%%%%%%%%%%%%%%%%%%%%%%%%%%%%%%%%%%%%%%%%%%%%%%%
%%%%%%%%%%%%%%%%%%%%%%%%%%%%%%%%%%%%%%%%%%%%%%%%%%%%%%%%%%%%%%%%%%%%%%%%%%%%%%%%%%%
\section{Introduction}
%%%%%%%%%%%%%%%%%%%%%%%%%%%%%%%%%%%%%%%%%%%%%%%%%%%%%%%%%%%%%%%%%%%%%%%%%%%%%%%%%%%
%%%%%%%%%%%%%%%%%%%%%%%%%%%%%%%%%%%%%%%%%%%%%%%%%%%%%%%%%%%%%%%%%%%%%%%%%%%%%%%%%%%
%%%%%%%%%%%%%%%%%%%%%%%%%%%%%%%%%%%%%%%%%%%%%%%%%%%%%%%%%%%%%%%%%%%%%%%%%%%%%%%%%%%

In this paper we consider the equation
\begin{equation}\label{bgk}
 \partial_t f^\eps +\partial_x (k(x) a(\xi) f^\eps)=\frac{\chi_{u^\eps}-f^\eps}{\eps},\quad
t>0, \, x \in\R, \,  \xi \in\R,
\end{equation}
with the initial condition 
\begin{equation} \label{ini}
f^\eps |_{t=0} = f_0, \mbox{ in }\R_x\times\R_\xi.
\end{equation}
Here $k$ is given by 
\begin{equation*}
k=k_L\1_{(-\infty,0)}+k_R\1_{(0,+\infty)},
\end{equation*}
where $\1_B$ is the characteristic function of a set $B$, $\xi\mapsto a(\xi)$ is a continuous function on $\R$ such that
\begin{equation}
\forall u\in [0,1],\int_0^u a(\xi)d\xi\geq 0,\quad \int_0^1 a(\xi)d\xi=0,
\label{HYPa}\end{equation}
and, in \refe{bgk}, $\chi_{u^\eps}$, the so-called {\it equilibrium function} associated to $f^\eps$ is defined by
\begin{equation*}
u^\eps(t,x)=\int_\R f^\eps(t,x,\xi)d\xi,\quad \chi_{\alpha}(\xi)=\1_{]0,\alpha[}(\xi)-\1_{]\alpha,0[}(\xi),
\end{equation*}
for $t>0,x\in\R,\xi\in\R,\alpha\in\R$.
\bigskip

Eq.~\refe{bgk} is the so-called BGK approximation to the scalar conservation law
\begin{equation}
\partial_t u+\partial_x(k(x)A(u))=0,\quad A(u)=\int_0^u a(\xi) d\xi.
\label{sclkdisc}\end{equation}

The flux $(x,u)\mapsto k(x)A(u)$ is discontinuous with respect to $x\in\R$, actually \refe{sclkdisc} is a prototype of scalar (first-order) conservation law with discontinuous flux function. In the last ten years, scalar conservation laws with discontinuous flux function have been extensively studied. We refer to the paper \cite{BurgerKarlsen08} for a comprehensive introduction to the subject and a complete list of references. Let us simply mention that the discontinuous character of the flux function gives rise to a multiplicity of weak solutions, even if traditional entropy conditions are imposed in the spatial domain apart from the discontinuity. An additional criterion has therefore to be given in order to select solutions in a unique way. For scalar conservation law under the general form $\partial_t u+\partial_x(B(x,u))=0$, where the function $B$ is discontinuous with respect to $x$, several criteria are possible~\cite{AdimurthiMishraGowda05}. For $B(x,u)=k(x)A(u)$ as above, the choice of entropy solution is unambiguous (see \cite{AdimurthiMishraGowda05}, Remark~4.4) and we consider here the criterion of selection first given in \cite{Towers01}. A kinetic formulation (in the spirit of \cite{LionsPerthameTadmor94}) equivalent to the entropy formulation in \cite{Towers01} has been given in \cite{BachmannVovelle06}. In particular, solutions given by this criterion are limits (a.e. and in $L^1$) of the solutions obtained by {\it monotone} regularization of the coefficient $k$ in \refe{sclkdisc}, {\it e.g.}
\begin{equation*}
k_\eps(x)=k_L\1_{x<-\eps}(x)+\left(\frac{k_R-k_L}{2\eps}x+\frac{k_R+k_L}{2}\right)\1_{-\eps\leq x\leq\eps}+k_R\1_{\eps<x},\quad \eps>0.
\end{equation*}

The kinetic formulation of scalar conservation laws is well adapted to the analysis of the (Perthame-Tadmor) BGK approximation of scalar conservation laws. Developed in~\cite{PerthameTadmor91}, this equation is a continuous version of the Transport-Collapse method of Brenier~\cite{Brenier81,Brenier83}.
BGK models have also been used for gas dynamics and the construction of numerical schemes.
See for example the book of Perthame \cite{Perthame02} for a survey of this field.
\medskip

Our purpose here is to apply the kinetic formulation of \cite{BachmannVovelle06} to show the convergence of the BGK approximation. 
To this aim, we first study the BGK equation in itself in Section~\ref{sec:bgk}. In Section~\ref{sec:sclkdisc}, we introduce the kinetic formulation for the limit problem. We also introduce a notion of generalized (kinetic) solution, Definition~\ref{def:gks}. We show that any generalized solution reduces to a mere solution, {\it i.e.} a solution in the sense of Def.~\ref{def:SOLSCL}. This theorem of ``reduction" is Theorem~\ref{th:redU}. Then in Section~\ref{sec:cvbgk}, we show that the BGK model converges to a generalized solution of \refe{sclkdisc} and, using Theorem~\ref{th:redU}, deduce the strong convergence of the BGK model to a solution of \refe{sclkdisc}, Theorem~\ref{th:cvbgk}.
\medskip

A key step of the whole proof of convergence is the result of reduction of Theorem~\ref{th:redU}. Its proof, given in Section~\ref{sec:proofredU}, is close to the proof of uniqueness of solutions given in \cite{BachmannVovelle06}. A minor difference is that we deal here with ge\-ne\-ra\-li\-zed solutions instead of ``kinetic process solutions". There is also a minor error in the proof given in \cite{BachmannVovelle06} (specifically, the remainder terms $R_{\alpha,\eps,\delta}$ and $Q_{\beta,\nu,\sigma}$ in Eq.~\refe{eq:ureg} and \refe{eq:vreg} of the present paper are missing in \cite{BachmannVovelle06}). We have therefore given a complete proof of Theorem~\ref{th:redU}.
\bigskip

We end this introduction with two remarks: 
\begin{itemize}
\item the BGK model provides an approximation of the entropy solutions to \refe{sclkdisc} by relaxation of the kinetic equation corresponding to \refe{sclkdisc}. A relaxation scheme of the Jin and Xin type applied directly to the original equation \refe{sclkdisc} has been developed in \cite{KarlsenKlingenbergRisebro04}.
\item in the last chapter of \cite{Bachmann05}, is derived the kinetic formulation of scalar conservation laws with discontinuous spatial dependence of the form $\partial_t u+\partial_x(B(x,u))=0$ (which are more general than \refe{sclkdisc}). We indicate (this would have to be proved rigorously), that in case where our approach {\it via} the BGK approximation was applied to this problem, the solutions obtained would be the type of entropy solutions considered in \cite{KarlsenRisebroTowers03}.
\end{itemize}

{\bf Notation} For $p,q\in[1,+\infty]$, 
we denote by $L^p_xL^q_\xi$ the space $L^p(\R_x;L^q(\R_\xi))$ and by $L^q_\xi L^p_x$ the space $L^q(\R_\xi;L^p(\R_x))$.
\smallskip

We also set $\sgn_+(s)=\1_{\{s>0\}}$, $\sgn_-(s)=-\1_{\{s\leq 0\}}$, $\sgn=\sgn_++\sgn_-$, $s\in\R$.
%%%%%%%%%%%%%%%%%%%%%%%%%%%%%%%%%%%%%%%%%%%%%%%%%%%%%%%%%%%%%%%%%%%%%%%%%%%%%%%%%%%
%%%%%%%%%%%%%%%%%%%%%%%%%%%%%%%%%%%%%%%%%%%%%%%%%%%%%%%%%%%%%%%%%%%%%%%%%%%%%%%%%%%
%%%%%%%%%%%%%%%%%%%%%%%%%%%%%%%%%%%%%%%%%%%%%%%%%%%%%%%%%%%%%%%%%%%%%%%%%%%%%%%%%%%
\section{The BGK equation}\label{sec:bgk}
%%%%%%%%%%%%%%%%%%%%%%%%%%%%%%%%%%%%%%%%%%%%%%%%%%%%%%%%%%%%%%%%%%%%%%%%%%%%%%%%%%%
%%%%%%%%%%%%%%%%%%%%%%%%%%%%%%%%%%%%%%%%%%%%%%%%%%%%%%%%%%%%%%%%%%%%%%%%%%%%%%%%%%%
%%%%%%%%%%%%%%%%%%%%%%%%%%%%%%%%%%%%%%%%%%%%%%%%%%%%%%%%%%%%%%%%%%%%%%%%%%%%%%%%%%%

%%%%%%%%%%%%%%%%%%%%%%%%%%%%%%%%%%%%%%%%%%%%%%%%%%%%%%%%%%%%%%%%%%%%%%%%%%%%%%%%%%%
\subsection{The balance equation}
%%%%%%%%%%%%%%%%%%%%%%%%%%%%%%%%%%%%%%%%%%%%%%%%%%%%%%%%%%%%%%%%%%%%%%%%%%%%%%%%%%%

By the change of variables $\tilde f^\eps(t,x,\xi)=e^{\frac{t}{\eps}}f^\eps(t,x,\xi)$, Eq.~\refe{bgk} rewrites as the balance equation
\begin{equation*}
\partial_t \tilde f^\eps+\partial_x(k(x)a(\xi)\tilde f^\eps)=\frac{e^{\frac{t}{\eps}}}{\eps}\chi_{u^\eps}
\end{equation*}
with (unknown dependent) source term $\frac{e^{\frac{t}{\eps}}}{\eps}\chi_{u^\eps}$. Hence, we first consider the following Cauchy Problem for the balance equation:
\begin{align}
\partial_t f+\partial_x(k(x)a(\xi)f)=&g,\quad t>0, x\in\R, \xi\in\R,\label{eq:balance}\\
f|_{t=0}= &f_0\quad\mbox{in }\R_x\times\R_\xi.\label{ic:balance}
\end{align}

\begin{proposition} Suppose that $k_R\cdot k_L>0$. Then Problem~\refe{eq:balance}-\refe{ic:balance} is well posed in $L^1_\xi L^p_x$, $1\leq p< +\infty$: for all $f_0\in L^1_\xi L^p_x $, $T>0$ and 
$g\in L^1(]0,T[;L^1_\xi L^p_x)$, there exists a unique $f\in C([0,T];L^1_\xi L^p_x)$ solving \refe{eq:balance} in $\mathcal{D}'(]0,T[\times\R_x\times\R_\xi)$ such that $f(0)=f_0$. Besides, we have
\begin{equation}
\|f(t)\|_{L^1_\xi L^p_x}\leq M_k\left(\|f_0\|_{L^1_\xi L^p_x}
+\int_0^t\|g(s)\|_{L^1_\xi L^p_x}ds\right),
\label{balanceLpL1}\end{equation}
where $M_k=\max\left(\frac{k_L}{k_R},\frac{k_R}{k_L}\right)$.
\label{prop:balance}\end{proposition}

{\bf Proof:} Since \refe{eq:balance} is linear, it is sufficient to solve the case $g=0$. The general case will follow from Duhamel's Formula. Assume without loss of generality $k_R,k_L>0$. Let $A_+:=\{\xi\in\R;a(\xi)>0\}$. Then, for fixed $\xi\in A_+$, and although $k$ is a discontinuous function, the O.D.E.
\begin{equation} 
\dot X(t,s,x,\xi)=k(X(t,s,x,\xi))a(\xi),\quad t\in\R,
\label{ODE}\end{equation}
with datum $X(s,s,x,\xi)=x$ has an obvious solution for $x\not=0$, given by $X(t,s,x,\xi)=x+(t-s)k_R a(\xi)$, $t>s$, when $x>0$, and by 
\begin{equation*}
X(t,s,x,\xi)=\left\{\begin{array}{l l l}
	x+(t-s)k_La(\xi) & \mbox{ if }&t<s+\frac{|x|}{k_La(\xi)},\\
	\frac{k_R}{k_L}x+(t-s)k_Ra(\xi) &\mbox{ if }& t>s +\frac{|x|}{k_La(\xi)},
\end{array}\right.
\end{equation*}
when $x<0$. Denoting by $s^+=\max(s,0)$, $s^-=s^+-s$ the positive and negative parts of $s\in\R$, and introducing
\begin{equation*}
\alpha_k(x)=\1_{\{x>0\}}+\frac{k_R}{k_L}\1_{\{x<0\}},
\end{equation*}
this can be summed up as 
\begin{equation}
X(t,s,x,\xi)=\{\alpha_k(x)x+(t-s)k_Ra(\xi)\}^+-\{x+(t-s)k_La(\xi)\}^-,\quad t>s.
\label{forward}\end{equation}
Similarly, we have, for the resolution of \refe{ODE} backward in time,
\begin{equation}
X(t,s,x,\xi)=\{x+(t-s)k_Ra(\xi)\}^+-\{\beta_k(x)x+(t-s)k_La(\xi)\}^-,\quad t<s,
\label{backward}\end{equation}
where
\begin{equation*}
\beta_k(x)=\frac{k_L}{k_R}\1_{\{x>0\}}+\1_{\{x<0\}}.
\end{equation*}
A similar computation in the case $a(\xi)\leq 0$ gives the solution to \refe{ODE} by \refe{forward} for $(t-s)a(\xi)\geq 0$, \refe{backward} for $(t-s)a(\xi)\leq 0$.
For the transport equation $(\partial_t+k(x)a(\xi)\partial_x)\varphi^*=0$, interpreted as 
\begin{equation*}
\frac{d\;}{dt}\varphi^*(t,X(t,s,x,\xi),\xi)=0,
\end{equation*}
this yields the solution
\begin{equation*}
\varphi^*(t,x,\xi)=\psi(X(T,t,x,\xi),\xi),
\end{equation*}
which satisfies the terminal condition $\varphi^*(T)=\psi$. We suppose in what follows that $\psi$ is independent on $\xi$, compactly supported and Lipschitz continuous. Then, a simple change of variable shows that, for every $t\in[0,T]$, for a.e. $\xi\in\R$, 
\begin{equation}
\|\varphi^*(t,\cdot,\xi)\|_{L^q_x}\leq M_k\|\psi\|_{L^q_x},\quad M_k=\max\left(\frac{k_L}{k_R},\frac{k_R}{k_L}\right), 1\leq q\leq +\infty.
\label{dualestimLq}\end{equation}
If $f\in C([0,T];L^1_\xi L^p_x)$ solves \refe{eq:balance}-\refe{ic:balance}, then, by duality (note that $\varphi^*$ is Lipschitz continuous and compactly supported in $x$ if $\psi$ is) we have, for $t\in[0,T]$, for a.e. $\xi\in\R$,
\begin{equation}
\int_\R f(T,x,\xi)\psi(x,\xi)dx=\int_\R f_0(x,\xi)\varphi^*(0,x,\xi)dx.
\label{duality}\end{equation}
In particular, the estimate \refe{dualestimLq} where $q=$ conjugate exponent of $p$ gives, for a.e. $\xi\in\R$,
$$\|f(T,\cdot,\xi)\|_{L^p_x}\leq M_k \|f_0(\cdot,\xi)\|_{L^p_x},$$
and then by Duhamel's principle, for $g \neq 0$,
\begin{equation}
\|f(T,\cdot,\xi)\|_{L^p_x}\leq M_k\left(\|f_0(\cdot,\xi)\|_{L^p_x}+\int_0^T \|g(t,\cdot,\xi)\|_{L^p_x}dt \right).
\label{balanceLp}\end{equation} 
The estimate~\refe{balanceLpL1} and uniqueness of the solution to \refe{eq:balance}-\refe{ic:balance} readily follows. Exis\-ten\-ce follows from \refe{forward}-\refe{backward}-\refe{duality}, from which one derives the explicit formula 
\begin{equation*}
f(t,x,\xi)=J(t,x,\xi)f_0(X(0,t,x,\xi),\xi),
\end{equation*}
the coefficient $J(t,x,\xi)$ being given by
\begin{equation*}
J(t,x,\xi)=\1_{\{x<0\}\cup\{x>tk_Ra(\xi)\}}+\frac{k_L}{k_R}\1_{\{0<x<tk_R a(\xi)\}}
\end{equation*}
if $a(\xi)>0$ and 
\begin{equation*}
J(t,x,\xi)=\1_{\{x<tk_La(\xi)\}\cup\{x>0\}}+\frac{k_R}{k_L}\1_{\{tk_L a(\xi)<x<0\}}
\end{equation*}
if $a(\xi)\leq 0$.

%%%%%%%%%%%%%%%%%%%%%%%%%%%%%%%%%%%%%%%%%%%%%%%%%%%%%%%%%%%%%%%%%%%%%%%%%%%%%%%%%%%
\subsection{The BGK equation}
%%%%%%%%%%%%%%%%%%%%%%%%%%%%%%%%%%%%%%%%%%%%%%%%%%%%%%%%%%%%%%%%%%%%%%%%%%%%%%%%%%%

Denote by ${\cal T}(t)f_0$ the solution to \refe{eq:balance}-\refe{ic:balance} with $g=0$, {\it i.e.}
\begin{equation*}
{\cal T}(t)f_0(x,\xi)=J(t,x,\xi)f_0(X(0,t,x,\xi),\xi),
\end{equation*}
$X$ given by \refe{forward}-\refe{backward}.

\begin{definition} Let $f_0\in L^1(\R_x\times\R_\xi)$, $T>0$. A function $f^\eps\in C([0,T];L^1(\R_x\times\R_\xi))$ is said to be a solution to \refe{bgk}-\refe{ini} if
\begin{equation}
f^\eps(t)=e^{-\frac{t}{\eps}}{\cal T}(t)f_0+
\frac{1}{\eps}\int_0^t e^{-\frac{s}{\eps}}{\cal T}(s)\chi_{u^\eps(t-s)}ds,\quad u^\eps=\int_\R f^\eps(\xi) d\xi,
\label{integralBGK}\end{equation}
for all $t\in[0,T]$.
\label{def:BGK}\end{definition}

\begin{theorem} Assume $k_R\cdot k_L>0$. Let $f_0\in L^1(\R_x\times\R_\xi)$, $T>0$. There exists a unique solution $f^\eps\in C([0,T];L^1(\R_x\times\R_\xi))$ to \refe{bgk}-\refe{ini}. Denoting by $S_\eps(t)f_0$ this solution, we have:
\begin{enumerate}
\item $\|(S_\eps(t)f_0^\natural-S_\eps(t)f_0^\flat)^+\|_{L^1(\R_x\times\R_\xi)}\leq M_k\|(f_0^\natural-f_0^\flat)^+\|_{L^1(\R_x\times\R_\xi)}$
\item $0\leq\sgn(\xi)f_0(x,\xi)\leq 1$ a.e. $\Rightarrow 0\leq\sgn(\xi)S_\eps(t)f_0(x,\xi)\leq 1$ a.e.
\item if $f_0=\chi_{u_0}$, $u_0\in L^\infty(\R)$, $0\leq u_0\leq 1$ a.e. then $0\leq S_\eps(t)f_0\leq\chi_1$.
\end{enumerate}
\label{th:BGK}\end{theorem}

{\bf Proof:} the change of variable $(t',x')=\eps(t,x)$ reduces \refe{bgk} to the same equation with $\eps=1$. We then have to solve $f=F(f)$ for
\begin{equation*}
F(f)(t):=e^{-t}{\cal T}(t)f_0+\int_0^t e^{-s}{\cal T}(s)\chi_{u(t-s)}ds,\quad u=\int_\R f(\xi) d\xi.
\end{equation*}
By \refe{balanceLpL1} and the identity
\begin{equation*}
\int_\R |\chi_u-\chi_v|(\xi)d\xi=|u-v|,\quad u,v\in\R,
\end{equation*}
we have $F\colon C([0,T];L^1_{x,\xi})\to C([0,T];L^1_{x,\xi})$ and $F$ is a $(1-e^{-T})$ contraction for the norm
\begin{equation*}
\|f\|=\sup_{t\in[0,T]}\|f(t)\|_{L^1(\R_x\times\R_\xi)}.
\end{equation*}
Indeed, we compute,
\begin{align*}
\|F(f^\natural)(t)-F(f^\flat)(t)\|_{L^1_{x,\xi}}\leq& \int_0^t e^{-s}\|{\cal T}(s)(\chi_{u^\natural (t-s)}-\chi_{u^\flat (t-s)})\|_{L^1_{x,\xi}}ds\\
=& \int_0^t e^{-s}\|\chi_{u^\natural (t-s)}-\chi_{u^\flat (t-s)}\|_{L^1_{x,\xi}}ds\\
=&\int_0^t e^{-s}\|u^\natural (t-s)-u^\flat (t-s)\|_{L^1_{x}}ds\\
\leq&\int_0^t e^{-s}\|f^\natural (t-s)-f^\flat (t-s)\|_{L^1_{x,\xi}}ds\\
\leq &\int_0^t e^{-s}ds \|f^\natural-f^\flat\|.
\end{align*}
By the Banach fixed point theorem, we obtain existence and uniqueness of the solution to \refe{bgk}-\refe{ini}.
Since $0\leq\sgn(\xi)\chi_u(\xi)\leq 1$ a.e. we have 
\begin{equation*}
0\leq\sgn(\xi)F(f)(t,x,\xi)\leq 1 \mbox{ a.e.}
\end{equation*} 
if $0\leq \sgn(\xi)f_0(x,\xi)\leq 1$ a.e. This proves the point {\it 2.} of the Theorem. The point {\it 1.} follows from the inequality
\begin{equation*}
\int_\R \sgn_+(f-g)(Q(f)-Q(g))d\xi\leq 0,\quad f,g\in L^1(\R_\xi),\quad Q(f):=\chi_{\int f d\xi}-f,
\end{equation*}
that is easy to check, and from the identity
\begin{equation*}
f(t)={\cal T}(t)f_0+\int_0^t {\cal T}(s)Q(f)(t-s)ds
\end{equation*}
for the solution to \refe{bgk}-\refe{ini}. If $f_0=\chi_{u_0}$, $0\leq u_0\leq 1$ a.e. then $0=\chi_0\leq f_0\leq\chi_1$. Hence the item {\it 3.} follows from {\it 1.} and the fact that any constant equilibrium function $\chi_\alpha$, $\alpha\in\R$ is solution to \refe{bgk}. \qed

%%%%%%%%%%%%%%%%%%%%%%%%%%%%%%%%%%%%%%%%%%%%%%%%%%%%%%%%%%%%%%%%%%%%%%%%%%%%%%%%%%%
%%%%%%%%%%%%%%%%%%%%%%%%%%%%%%%%%%%%%%%%%%%%%%%%%%%%%%%%%%%%%%%%%%%%%%%%%%%%%%%%%%%
%%%%%%%%%%%%%%%%%%%%%%%%%%%%%%%%%%%%%%%%%%%%%%%%%%%%%%%%%%%%%%%%%%%%%%%%%%%%%%%%%%%
\section{The limit problem}\label{sec:sclkdisc}
%%%%%%%%%%%%%%%%%%%%%%%%%%%%%%%%%%%%%%%%%%%%%%%%%%%%%%%%%%%%%%%%%%%%%%%%%%%%%%%%%%%
%%%%%%%%%%%%%%%%%%%%%%%%%%%%%%%%%%%%%%%%%%%%%%%%%%%%%%%%%%%%%%%%%%%%%%%%%%%%%%%%%%%
%%%%%%%%%%%%%%%%%%%%%%%%%%%%%%%%%%%%%%%%%%%%%%%%%%%%%%%%%%%%%%%%%%%%%%%%%%%%%%%%%%%

Assume $f_0=\chi_{u_0}$ with $u_0\in L^\infty(\R)$, $0\leq u_0\leq 1$ a.e. Set 
\begin{equation} \label{hypsurA}
A(u)=\int_0^u a(\xi)\1_{[0,1]}(\xi) d\xi.
\end{equation} 
Note that by \refe{HYPa}, we have $A\geq 0$ and $A$ vanishes outside the interval $[0,1]$. We expect the solution $f^\eps$ to \refe{bgk}-\refe{ini} to converge to the solution $u$ of the first-order scalar conservation law
\begin{equation}
\partial_t u+\partial_x(k(x)A(u))=0,\quad t>0, x\in\R,
\label{SCL}\end{equation}
with initial datum
\begin{equation}
u(0,x)=u_0(x),\quad x\in\R.
\label{IC}\end{equation}

For a fixed $T>0$, set $Q=]0,T[\times\R_x$.
\medskip

%Def sol
\begin{definition}[Solution] Let $u_{0}\in L^\infty(\R)$, $0\leq u_0\leq 1$ a.e. A function $u\in L^\infty(Q)$ is said to be a (kinetic) solution to~\refe{SCL}-\refe{IC} if there exists non-negative measures $m_\pm$ on $[0,T]\times\R\times\R$ such that
\begin{itemize}
\item $m_+$ is supported in $[0,T]\times\R\times]-\infty,1]$, $m_-$ is supported in $[0,T]\times\R\times[0,+\infty[$,
\item for all $\psi\in C^\infty_c([0,T[\times\R\times\R)$,
\begin{multline}
\int_{Q}\int_\R h_\pm (\partial_t\psi+k(x) a(\xi)\partial_x\psi) d\xi dt dx\\
+\int_{\R}\int_\R h_{0,\pm}\psi(0,x,\xi) d\xi dx -(k_L-k_R)^\pm\int_0^T\int_\R a(\xi) \psi(t,0,\xi) d\xi dt\\
=\int_{Q}\int_\R\partial_\xi\psi dm_\pm(t,x,\xi)
\label{eq:kineticu}\end{multline}
where $h_\pm(t,x,\xi)=\sgn_\pm(u(t,x)-\xi)$, $h_{0,\pm}(x,\xi)=\sgn_\pm(u_0(x)-\xi)$.
\end{itemize}
\label{def:SOLSCL}\end{definition}

%Bound Sol
\begin{proposition}[Bound in $L^\infty$] Let $u_{0}\in L^\infty(\R)$, $0\leq u_0\leq 1$ a.e. If $u\in L^\infty(Q)$ is a kinetic solution to~\refe{SCL}-\refe{IC}, then $0\leq u\leq 1$ a.e.
\label{prop:boundu}\end{proposition}

{\bf Proof:} Consider the kinetic formulation~\refe{eq:kineticu} for $h_+$ with a test function 
\begin{equation*}
\psi(t,x,\xi)=\varphi(t,x)\mu(\xi).
\end{equation*} 
If $\mu$ is supported in $]1,+\infty[$, two terms cancel:
\begin{equation*}
\int_{\R}\int_\R h_{0,+}\psi(0,x,\xi) d\xi dx=
\int_\R\int_\R\1_{1 \geq u_0(x)>\xi}\varphi(0,x) \1_{\xi >1}\mu(\xi)d\xi dx=0
\end{equation*}
and 
\begin{equation*}
\int_{Q}\int_\R\partial_\xi\psi dm_+(t,x,\xi)=0
\end{equation*}
by the hypothesis on the support of $m_+$. Hence we have
\begin{multline*}
\int_{Q}\int_\R h_+ (\partial_t\varphi+k(x) a(\xi)\partial_x\varphi)\mu(\xi) d\xi dt dx\\
-(k_L-k_R)^+\int_0^T\int_\R a(\xi) \varphi(t,0)\mu(\xi) d\xi dt=0.
\end{multline*}
A step of approximation and regularization shows that we can take $\mu(\xi)=\1_{\xi>1}$ in this equation. Since 
\begin{equation*}
\int_1^{+\infty}a(\xi)d\xi=A(+\infty)-A(1)=0-0=0,
\end{equation*}
and
\begin{equation*}
\int_1^{+\infty} h_+(t,x,\xi) d\xi=\int_1^{+\infty} \1_{\xi <u(t,x)} d\xi =(u(t,x)-1)^+,
\end{equation*}
\begin{multline*}
\int_1^{+\infty} h_+(t,x,\xi) a(\xi) d\xi
=\int_1^{+\infty} \1_{\xi <u(t,x)} a(\xi) d\xi \\ =
\sgn_+(u(t,x)-1) \int_1^{u(t,x)}a(\xi) d\xi =
\sgn_+(u(t,x)-1) (A(u(t,x)-A(1)) ,
\end{multline*}
we obtain
\begin{equation*}
\int_{Q} (u-1)^+\partial_t\varphi+k(x) \sgn_+(u-1)(A(u)-A(1))\partial_x\varphi dt dx =0.
\end{equation*}
It is then classical to deduce that $(u-1)^+=0$ a.e. (see the end of the proof of Proposition~\ref{prop:compfg}, after \refe{U5}), {\it i.e.} $u\leq 1$ a.e. Similarly, we show $u\geq 0$ a.e. \qed
\medskip

Our aim is to prove the uniqueness of the solution to \refe{SCL}-\refe{IC}. Actually, more than mere uniqueness of the solution to \refe{SCL}-\refe{IC}, we will show a result of reduction/uniqueness (see Theorem~\ref{th:redU}) of generalized kinetic solution. To this purpose, let us recall that a Young measure $Q\to\R$ is a measurable mapping $(t,x)\mapsto\nu_{t,x}$ from $Q$ into the space of probability (Borel) measures on $\R$. The mapping is measurable in the sense that for each Borel subset $A$ of $\R$, $(t,x)\mapsto\nu_{t,x}(A)$ is measurable $Q\to\R$. Let us also introduce the following notation: if $f\in L^1(Q\times\R)$, we set 
\begin{equation*}
f_\pm(y,\xi)=f(y,\xi)-\sgn_\mp(\xi),\quad y\in Q,\xi\in\R.
\end{equation*}
This is consistent with the notations used in Def.~\ref{def:SOLSCL} in the case $f=\chi_u$.

%Def gks
\begin{definition}[Generalized solution] Let $u_{0}\in L^\infty(\R)$, $0\leq u_0\leq 1$ a.e. A function $f\in L^1(Q\times\R_\xi)$ is said to be a generalized (kinetic) solution to~\refe{SCL}-\refe{IC} if 
\begin{equation*}
0\leq f\leq\chi_1\mbox{ a.e., }-\partial_\xi f_+\mbox{ is a Young measure }Q\to\R,
\end{equation*}
and if there exists non-negative measures $m_\pm$ on $[0,T]\times\R\times\R$ such that
\begin{itemize}
\item $m_+$ is supported in $[0,T]\times\R\times]-\infty,1]$, $m_-$ is supported in $[0,T]\times\R\times[0,+\infty[$,
\item for all $\psi\in C^\infty_c([0,T[\times\R\times\R)$,
\begin{multline}
\int_{Q}\int_\R f_\pm (\partial_t\psi+k(x) a(\xi)\partial_x\psi) d\xi dt dx\\
+\int_{\R}\int_\R f_{0,\pm}\psi(0,x,\xi) d\xi dx-(k_L-k_R)^\pm\int_0^T\int_\R a(\xi) \psi(t,0,\xi) d\xi dt\\
=\int_{Q}\int_\R\partial_\xi\psi dm_\pm(t,x,\xi)
\label{eq:kineticf}\end{multline}
where $f_{0,\pm}(x,\xi)=\sgn_\pm(u_0(x)-\xi)$.
\end{itemize}
\label{def:gks}\end{definition}

\begin{theorem}[Reduction, Uniqueness] Let $u_{0}\in L^\infty(\R)$, $0\leq u_0\leq 1$ a.e. Problem~\refe{SCL}-\refe{IC} admits at most one solution. Besides, any generalized solution is actually a solution: if $f\in L^1(Q\times\R_\xi)$ is a generalized solution to \refe{SCL}-\refe{IC}, then there exists $u\in L^\infty(Q)$ such that $f=\chi_u$.
\label{th:redU}\end{theorem}

To prepare the proof of Theorem~\ref{th:redU}, we first have to analyze the formulation \refe{eq:kineticf} and the behavior of $f$ at $t=0$ and $x=0$.

%%%%%%%%%%%%%%%%%%%%%%%%%%%%%%%%%%%%%%%%%%%%%%%%%%%%%%%%%%%%%%%%%%%%%%%%%%%%%%%%%%
\subsection{Weak traces}
%%%%%%%%%%%%%%%%%%%%%%%%%%%%%%%%%%%%%%%%%%%%%%%%%%%%%%%%%%%%%%%%%%%%%%%%%%%%%%%%%%

Introduce the cut-off function
\begin{equation}
\omega_\eps(s)=\int_0^{|s|}\rho_\eps(r)dr,\quad \rho_\eps(s)=\eps^{-1}\rho(\eps^{-1}s),\quad s\in\R,
\label{defrho}\end{equation}
where $\rho\in C^\infty_c(\R)$ is a non-negative function with total mass $1$ compactly supported in $]0,1[$. We have the following proposition.

%weak traces
\begin{proposition}[Weak traces] Let $f\in L^\infty(Q\times\R_\xi)$ be a generalized solution to \refe{SCL}-\refe{IC}. There exists $f^{\tau_0}_\pm\in L^2(\R\times\R)$, $F_\pm\in L^2(]0,T[\times\R)$ and a sequence $(\eta_n)\downarrow 0$ such that, for all $\varphi\in L^2_c(\R\times\R)$, for all $\theta\in L^2_c(]0,T[\times\R)$ (the subscript $c$ denotes compact support),
\begin{align}
\int_Q\int_\R f_\pm(t,x,\xi)\omega_{\eta_n}'(t)\varphi(x,\xi)d\xi dt dx\to&\int_\R\int_\R f^{\tau_0}_\pm(x,\xi)\varphi(x,\xi) d\xi dx,\label{tracetime}\\
\int_Q\int_\R f_\pm(t,x,\xi)k(x) a(\xi)\omega_{\eta_n}'(x)\theta(t,\xi)d\xi dt dx\to&\int_0^T\int_\R F_\pm(t,\xi)\theta(t,\xi) d\xi dt\label{tracespace}
\end{align}
as $n\to+\infty$. Besides, there exists non-negative measures $m^{\tau_0}_\pm$, $\bar m_\pm$ on $\R^2$ and $[0,T]\times\R$ respectively such that:
\begin{itemize}
\item $m^{\tau_0}_+$ ({\it resp.} $\bar m_+$) is supported in $\R\times]-\infty,1]$ ({\it resp.} $[0,T]\times]-\infty,1]$), $m^{\tau_0}_-$ ({\it resp.} $\bar m_-$) is supported in $\R\times[0,+\infty[$ ({\it resp.} $[0,T]\times[0,+\infty[$),
\item for all $\varphi\in C^\infty_c(\R^2)$, $\theta\in C^\infty_c([0,T[\times\R)$,
\begin{align}
\int_{\R^2}f^{\tau_0}_\pm \varphi dx d\xi=&\int_{\R^2}f_{0,\pm} \varphi dx d\xi-\int_{\R^2}\partial_\xi\varphi dm^{\tau_0}_\pm(x,\xi),\label{meastracetime}\\
\int_0^T\int_{\R}F_\pm \theta  d\xi dt=&-(k_L-k_R)^\pm \int_0^T\int_{\R}a(\xi) \theta
 d\xi dt\nonumber\\
&-\int_0^T\int_{\R}\partial_\xi\theta d\bar m_\pm(t,\xi).\label{meastracespace}
\end{align}
\end{itemize}
\label{prop:trace}\end{proposition}

{\bf Proof: } The first part of the proposition does not use the fact that $f$ is solution. Indeed, since $|f_\pm|\leq 2$, we have 
\begin{equation*}
\left|\int_0^Tf_\pm(t,x,\xi)\omega_{\eta}'(t) dt \right|\leq 2\int_0^T |\omega_\eta'(t)|dt=2\int_0^T\rho_\eta(t)dt\leq 2,
\end{equation*}
for all $(x,\xi)\in\R^2$. This gives in particular a bound in $L^2(K)$, $K$ compact of $\R^2$ on $\int_0^T f_\pm(t,\cdot)\omega_\eta'(t)dt$, hence existence of a subsequence that converges weakly in $L^2(K)$. Writing $\R^2$ as an increasing countable union of compact sets and using a diagonal process, we obtain \refe{tracetime}. The proof of \refe{tracespace} is similar.
To obtain \refe{meastracetime}, apply the formulation \refe{eq:kineticf} to $\psi(t,x,\xi)=\varphi(x,\xi)(1-\omega_{\eta_n}(t))$. We obtain \refe{meastracetime} by using \refe{tracetime} and setting
\begin{equation*}
\int_{\R^2}\varphi dm^{\tau_0}_\pm(x,\xi)=\lim_{n\to+\infty}\int_Q\int_{\R} \varphi(x,\xi)(1-\omega_{\eta_n}(t)) dm_\pm(t,x,\xi)
\end{equation*}
for all non-negative $\varphi\in C_c(\R^2)$: the limit is well defined since the argument is monotone in $n$ and it defines a non-negative functional on $C_c(\R^2)$ which is re\-pre\-sen\-ted by a non-negative Radon measure. Similarly, applying the  formulation \refe{eq:kineticf} to $\psi(t,x,\xi)=\theta(t,\xi)(1-\omega_{\eta_n}(x))$, we obtain \refe{meastracespace} with
\begin{equation*}
\int_0^T\int_{\R}\theta d\bar m_\pm(t,\xi)=\lim_{n\to+\infty}\int_{Q}\int_\R \theta(t,\xi)(1-\omega_{\eta_n}(x)) dm_\pm(t,x,\xi)
\end{equation*}
for all non-negative $\theta\in C_c([0,T]\times\R)$. \qed
\medskip

{\bf Remark:} Since $0\leq f\leq \chi_1$, \refe{tracetime} shows that $f^{\tau_0}_+$, {\it resp.} $f^{\tau_0}_-$, is supported in
$\R\times]-\infty,1]$, {\it resp.} $\R\times[0,+\infty[$. Similarly, $F_+$, {\it resp.} $F_-$, is supported in
$[0,T]\times]-\infty,1]$, {\it resp.} $[0,T]\times[0,+\infty[$. We use this remark to show the following 

%corollary
\begin{corollary} For all $\varphi_-\in L^\infty(\R^2)$ supported in $[-R,R]\times [-R,+\infty[$ ($R>0$) such that $\partial_\xi\varphi_-\leq 0$ (in the sense of distributions), we have 
\begin{equation}
\lim_{n\to+\infty}\int_Q\int_\R f_+\omega_{\eta_n}'(t)\varphi_-(x,\xi)d\xi dt dx\geq \int_{\R^2}f_{0,+} \varphi_- dx d\xi.
\label{tracetimeext}\end{equation}
For all $\theta_-\in L^\infty(]0,T[\times\R)$ supported in $[0,T]\times [-R,+\infty[$ ($R>0$) such that $\partial_\xi\theta_-\leq 0$ (in the sense of distributions), we have 
\begin{equation}
\lim_{n\to+\infty} \int_Q\int_\R f_+k(x) a(\xi)\omega_{\eta_n}'(x)\theta_-(t,\xi)d\xi dt dx\geq -(k_L-k_R)^+ \int_0^T\int_{\R}a(\xi) \theta_- d\xi dt.
\label{tracespaceext}\end{equation}
\end{corollary}

{\bf Proof: } Note first that each term in \refe{tracetimeext} is well defined by the remark above and that, by \refe{tracetime},
\begin{equation*}
\lim_{n\to+\infty}\int_Q\int_\R f_+(t,x,\xi)\omega_{\eta_n}'(t)\varphi_-(x,\xi)d\xi dt dx=\int_\R\int_\R f_+^{\tau_0}\varphi_-  d\xi dx.
\end{equation*}
By regularization (parameter $\eps$) and truncation (parameter $M$), we have
\begin{equation*}
\int_{\R^2}(f^{\tau_0}_+ -f_{0,+}) \varphi_- dx d\xi=\int_{\R^2}(f^{\tau_0}_+ -f_{0,+}) \varphi_-^{\eps,M} dx d\xi+\eta(\eps,M),
\end{equation*}
where $\limd_{\eps\to0,M\to+\infty}\eta(\eps,M)=0$. More precisely, we set
\begin{equation*}
\varphi_-^{\eps,M}=(\varphi_-*\psi_\eps)\times\chi_M,
\end{equation*}
where $\psi_\eps$ is a (smooth, compactly supported) approximation of the unit on $\R^2$ and $\chi_M$ is a smooth, non-increasing function such that $\chi_M\equiv 1$ on $]-\infty,M]$, $\chi_M\equiv 0$ on $[M+1,+\infty[$. Apply \refe{meastracetime} to $\varphi_-^{\eps,M}$ to obtain
\begin{equation*}
\int_{\R^2}(f^{\tau_0}_+ -f_{0,+}) \varphi_- dx d\xi=-\int_{\R^2}\partial_\xi\varphi_-^{\eps,M} dm^{\tau_0}_+(x,\xi)+\eta(\eps,M).
\end{equation*}
For $M>R+1$ and $\eps<1$, we have $\varphi_-^{\eps,M}=\varphi_-*\psi_\eps$, hence $\partial_\xi \varphi_-^{\eps,M}\leq 0$. It follows that
\begin{equation*}
\int_{\R^2}(f^{\tau_0}_+ -f_{0,+}) \varphi_- dx d\xi\geq\eta(\eps,M),
\end{equation*}
for $M>R+1$, $\eps<1$. At the limit $M\to+\infty$, $\eps\to 0$, we obtain \refe{tracetimeext}. The proof of \refe{tracespaceext} is similar. \qed

%%%%%%%%%%%%%%%%%%%%%%%%%%%%%%%%%%%%%%%%%%%%%%%%%%%%%%%%%%%%%%%%%%%%%%%%%%%%%%%%%%
\subsection{Proof of Theorem~\ref{th:redU}}\label{sec:proofredU}
%%%%%%%%%%%%%%%%%%%%%%%%%%%%%%%%%%%%%%%%%%%%%%%%%%%%%%%%%%%%%%%%%%%%%%%%%%%%%%%%%%

Our aim is to show the following 
\begin{proposition} Let $u_0,v_0\in L^\infty(\R)$, $0\leq u_0,v_0\leq 1$ a.e. and let $f$, {\it resp} $g$, be a generalized solution to \refe{SCL}-\refe{IC} with datum $u_0$, {\it resp.} $v_0$. Let $M=\supd_{x\in\R,\xi\in[0,1]}|k(x)a(\xi)|$. Then we have, for $R>0$,
\begin{equation}
\frac{1}{T}\int_0^T\int_{\{|x|<R\}}\int_\R -f_+g_- d\xi dx dt \leq \int_{\{|x|<R+MT\}}(u_0-v_0)^+ dx.  
\label{eq:compfg}\end{equation}
\label{prop:compfg}\end{proposition}

{\bf Remark:} 
In case $f=\chi_u$, $g=\chi_v$, we have $\int_\R -f_+ g_-d\xi=(u-v)^+$, hence \refe{eq:compfg} gives uniqueness of the solution to \refe{SCL}-\refe{IC} (more precisely, it gives the $L^1$-contraction with averaging in time and the comparison result $u_0\leq v_0$ a.e. $\Rightarrow u\leq v$ a.e.).

{\bf Remark:} 
To obtain the second part of Theorem~\ref{th:redU}, we apply \refe{eq:compfg} with $g=f$ to obtain
\begin{equation}
\int_0^T\int_{\{|x|<R\}}\int_\R -f_+f_- d\xi dx dt\leq 0.
\label{eq:compfg2}\end{equation}
Since $0\leq f\leq\chi_1$, we have $f_+\geq 0$ a.e. and $f_-\leq 0$ a.e. We deduce from \refe{eq:compfg2} that $f_+f_-=0$ a.e. Let $\nu_{t,x}$ denote the Young measure $-\partial_\xi f_+$: we have $\partial_\xi f_-=\partial_\xi f -\delta_0=\partial_\xi f_+$ and, by examination of the values at $\xi=\pm\infty$ of $f_\pm$, for a.e. $(t,x)\in Q$,
\begin{equation*}
f_+(t,x,\xi)=\nu_{t,x}(\xi,+\infty),\quad f_-(t,x,\xi)=-\nu_{t,x}(-\infty,\xi).
\end{equation*}
But then, the relation $f_+f_-=0$ implies that $\nu_{t,x}$ is a Dirac mass at, say, $u(t,x)$. By measurability of $\nu$, $u$ is measurable and $f=\chi_u$.
\bigskip

{\bf Proof of Proposition~\ref{prop:compfg}:} Since $f_+$ and $g_-$ satisfy
\begin{multline}
\int_{Q}\int_\R f_+ (\partial_t\psi+k(x) a(\xi)\partial_x\psi) d\xi dt dx\\
+\int_{\R}\int_\R f_{0,+}\psi(0,x,\xi) d\xi dx-(k_L-k_R)^+\int_0^T\int_\R a(\xi) \psi(t,0,\xi)  d\xi dt\\
=\int_{Q}\int_\R\partial_\xi\psi dm_+(t,x,\xi)
\label{eq:kineticuu}\end{multline}
and
\begin{multline}
\int_{Q}\int_\R g_- (\partial_t\psi+k(x) a(\xi)\partial_x\psi) d\xi dt dx\\
+\int_{\R}\int_\R g_{0,-}\psi(0,x,\xi) d\xi dx-(k_L-k_R)^-\int_0^T\int_\R a(\xi) \psi(t,0,\xi) d\xi dt\\
=\int_{Q}\int_\R\partial_\xi\psi dp_-(t,x,\xi)
\label{eq:kineticvv}\end{multline}
for all $\psi\in C^\infty_c([0,T[\times\R\times\R)$ (here $g_{0,-}=\sgn_-(v_0-\xi)$ and $p_-$ is a non-negative measure on $[0,T]\times\R\times\R$ supported in $[0,T]\times\R\times[0,+\infty[$), it is possible to obtain an estimate for $-f_+g_-$ by setting $\psi=-g_-\varphi$ in \refe{eq:kineticuu} and $\psi=f_+\varphi$ in \refe{eq:kineticvv} ($\varphi$ being a given test function) and adding the result. This requires first, however, a step of regularization.
\medskip

%Regularization
{\bf Step 1. Regularization.} Let $\rho_{\alpha,\eps,\delta}$ denote the approximation of the unit on $\R^3$ given by
\begin{equation*}
\rho_{\alpha,\eps,\delta}(t,x,\xi)=\rho_\alpha(t)\rho_\eps(x)\rho_\delta(\xi),\quad (t,x,\xi)\in\R^3,
\end{equation*}
where $\rho_\eps$ is defined in \refe{defrho}. Let $\psi\in C^\infty_c([0,T[\times\R\times\R)$ be compactly supported in $]0,T[\times\R\setminus\{0\}\times\R$. 
Use $\psi*\rho_{\alpha,\eps,\delta}$ as a test function in \refe{eq:kineticuu} and Fubini's theorem to obtain
\begin{multline*}
\int_{Q}\int_\R f_+^{\alpha,\eps,\delta} (\partial_t\psi+k(x) a(\xi)\partial_x\psi) d\xi dt dx\\
+\int_{\R}\int_\R f_{0,+}\psi*\rho_{\alpha,\eps,\delta}(0,x,\xi) d\xi dx-(k_L-k_R)^+\int_0^T\int_\R a(\xi) \psi*\rho_{\alpha,\eps,\delta}(t,0,\xi) d\xi dt\\
=\int_{Q}\int_\R\partial_\xi\psi dm_+^{\alpha,\eps,\delta}(t,x,\xi)+R_{\alpha,\eps,\delta}(\psi),
\end{multline*}
where $f_+^{\alpha,\eps,\delta}:=f_+*\check{\rho}_{\alpha,\eps,\delta}$, $m_+^{\alpha,\eps,\delta}:=m_+*\check{\rho}_{\alpha,\eps,\delta}$ and
\begin{equation*}
R_{\alpha,\eps,\delta}(\psi)=\int_Q\int_\R f_+[k(x) a(\xi)(\partial_x\psi)*\rho_{\alpha,\eps,\delta}-(k(x) a(\xi)\partial_x\psi)*\rho_{\alpha,\eps,\delta}]d\xi dt dx.
\end{equation*}
Here we have denoted $\check{\rho}(t,x,\xi)=\rho(-t,-x,-\xi)$. Also observe that, implicitly, we have extended $f_+$ by $0$ outside $[0,T]$ since, {\it e.g.}
\begin{align*}
\int_0^T f_+(t) \psi*\rho_\alpha(t) dt=& \int_0^T\int_\R f_+(t)\psi(s)\rho_\alpha(t-s)ds dt\\
=& \int_\R \psi(s) \int_0^T f_+(t)\check{\rho}_\alpha(s-t) dt ds.
\end{align*}
Since $\psi$ is supported in $]0,T[\times\R\setminus\{0\}\times\R$, we have, for $\alpha,\eps$ small enough,
\begin{align*}
\int_{\R}\int_\R f_{0,+}\psi*\rho_{\alpha,\eps,\delta}(0,x,\xi) d\xi dx=&0,\\
\int_0^T\int_\R a(\xi) \psi*\rho_{\alpha,\eps,\delta}(t,0,\xi)  d\xi dt=&0,\\
\end{align*}
and
\begin{equation*}
R_{\alpha,\eps,\delta}(\psi)=\int_Q\int_\R f_+ k(x)[a(\xi)(\partial_x\psi)*\rho_{\alpha,\eps,\delta}-(a(\xi)\partial_x\psi)*\rho_{\alpha,\eps,\delta}]d\xi dt dx.
\end{equation*}
We deduce 
\begin{multline}
\int_{Q}\int_\R f_+^{\alpha,\eps,\delta} (\partial_t\psi+k(x) a(\xi)\partial_x\psi) d\xi dt dx\\
=\int_{Q}\int_\R\partial_\xi\psi dm_+^{\alpha,\eps,\delta}(t,x,\xi)+R_{\alpha,\eps,\delta}(\psi).
\label{eq:ureg}\end{multline}
A similar work on $g_-$ gives 
\begin{multline}
\int_{Q}\int_\R g_-^{\beta,\nu,\sigma} (\partial_t\psi+k(x) a(\xi)\partial_x\psi) d\xi dt dx\\
=\int_{Q}\int_\R\partial_\xi\psi dp_-^{\beta,\nu,\sigma}(t,x,\xi)+Q_{\beta,\nu,\sigma}(\psi),
\label{eq:vreg}\end{multline}
where
\begin{equation*}
Q_{\beta,\nu,\sigma}(\psi)=\int_Q\int_\R g_- k(x)[a(\xi)(\partial_x\psi)*\rho_{\beta,\nu,\sigma}-(a(\xi)\partial_x\psi)*\rho_{\beta,\nu,\sigma}]d\xi dt dx.
\end{equation*}

%Mix
{\bf Step 2. Equation for $-f_+^{\alpha,\eps,\delta}g_-^{\beta,\nu,\sigma}$.} Let $\varphi\in C^\infty_c([0,T[\times\R)$ be non-negative and compactly supported in $]0,T[\times \R\setminus\{0\}$. Notice that $\varphi$ does not depend on $\xi$. Set $\psi=-\varphi g_-^{\beta,\nu,\sigma}$ in \refe{eq:ureg}, $\psi=-\varphi f_+^{\alpha,\eps,\delta}$ in \refe{eq:vreg}. Since
\begin{equation*}
f\partial_t(\varphi g)+g\partial_t(\varphi f)=fg\partial_t\varphi+\partial_t(\varphi fg),
\end{equation*} 
we obtain by addition of the resulting equations
\begin{multline*}
\int_{Q}\int_\R -f_+^{\alpha,\eps,\delta}g_-^{\beta,\nu,\sigma} (\partial_t\varphi+k(x) a(\xi)\partial_x\varphi) d\xi dt dx\\
=-\int_{Q}\varphi \int_\R\partial_\xi f_+^{\alpha,\eps,\delta}dp_-^{\beta,\nu,\sigma}(t,x,\xi)
+\partial_\xi g_-^{\beta,\nu,\sigma}dm_+^{\alpha,\eps,\delta}(t,x,\xi)\\
+R_{\alpha,\eps,\delta}(-\varphi g_-^{\beta,\nu,\sigma})+Q_{\beta,\nu,\sigma}(-\varphi f_+^{\alpha,\eps,\delta}).
\end{multline*}
Notice that the term
\begin{equation*}
-\int_{Q}\varphi \int_\R\partial_\xi f_+^{\alpha,\eps,\delta}dp_-^{\beta,\nu,\sigma}(t,x,\xi)
+\partial_\xi g_-^{\beta,\nu,\sigma}dm_+^{\alpha,\eps,\delta}(t,x,\xi)
\end{equation*}
is well defined since the intersection of the supports of the functions $f_+^{\alpha,\eps,\delta}$ and $p_-^{\beta,\nu,\sigma}$ ({\it resp.} $f_-^{\beta,\nu,\sigma}$ and $m_+^{\alpha,\eps,\delta}$) is compact. Actually, this term is non-negative since
$p_-^{\beta,\nu,\sigma},m_+^{\alpha,\eps,\delta}\geq 0$ and 
$\partial_\xi f_+^{\alpha,\eps,\delta},\partial_\xi g_-^{\beta,\nu,\sigma}\leq 0$. We thus have
\begin{multline}
\int_{Q}\int_\R -f_+^{\alpha,\eps,\delta}g_-^{\beta,\nu,\sigma} (\partial_t\varphi+k(x) a(\xi)\partial_x\varphi) d\xi dt dx\\
\geq R_{\alpha,\eps,\delta}(-\varphi g_-^{\beta,\nu,\sigma})+Q_{\beta,\nu,\sigma}(-\varphi f_+^{\alpha,\eps,\delta}).
\label{U1}\end{multline}
It is easily checked that 
\begin{equation*}
R_{\alpha,\eps,\delta}(-\varphi j_-^{\beta,\nu,\sigma})=\mathcal{O}(\nu^{-1}\delta),\quad Q_{\beta,\nu,\sigma}(-\varphi h_+^{\alpha,\eps,\delta})=\mathcal{O}(\eps^{-1}\sigma),
\end{equation*} 
hence
\begin{equation*}
\lim_{\delta,\sigma\to 0}R_{\alpha,\eps,\delta}(-\varphi g_-^{\beta,\nu,\sigma})+Q_{\beta,\nu,\sigma}(-\varphi f_+^{\alpha,\eps,\delta})=0.
\end{equation*}
At the limit $\delta,\sigma\to 0$ in \refe{U1}, we conclude that
\begin{equation}
\int_{Q}\int_\R -f_+^{\alpha,\eps}g_-^{\beta,\nu} (\partial_t\varphi+k(x) a(\xi)\partial_x\varphi) d\xi dt dx\geq 0.
\label{U2}\end{equation}

%Traces
{\bf Step 3. Traces.} Suppose that $k_L<k_R$. We then pass to the limit $\eps,\alpha\to 0$ in \refe{U2} to obtain
\begin{equation}
\int_{Q}\int_\R -f_+g_-^{\beta,\nu} (\partial_t\varphi+k(x) a(\xi)\partial_x\varphi) d\xi dt dx\geq 0.
\label{U3}\end{equation}
{\it Note} that in the opposite case $k_L>k_R$, and with our method of proof, we would {\it first} pass to the limit on $\beta,\nu$. Let us now remove the hypothesis that $\varphi$ vanishes at $t=0$: suppose that $\psi\in C^\infty_c([0,T[\times\R)$ is non-negative and supported in $[0,T[\times\R\setminus\{0\}$ and apply \refe{U3} to $\varphi(t,x)=\psi(t,x)\omega_{\eta_n}(t)$. We have
\begin{multline} 
\int_{Q}\int_\R -f_+g_-^{\beta,\nu} \omega_{\eta_n}(t)(\partial_t\psi
+k(x) a(\xi)\partial_x\psi) d\xi dt dx\\
+\int_Q\int_\R -f_+g_-^{\beta,\nu}\psi(t,x)\omega_{\eta_n}'(t) d\xi dt dx  \geq 0. \label{eqi}
\end{multline}
By \refe{tracetimeext} applied with $\varphi_-(x,\xi)=g_-^{\beta,\nu}(0,x,\xi)\psi(0,x)$, we obtain
\begin{multline*}
\lim_{n\to+\infty}\int_Q\int_\R f_+g_-^{\beta,\nu}(0,x,\xi)\psi(0,x)\omega_{\eta_n}'(t) d\xi dt dx\\
\geq \int_\R\int_\R f_{0,+}g_-^{\beta,\nu}(0,x,\xi)\psi(0,x) d\xi dx.
\end{multline*}
Now 
$f_+(t,x,\xi)g_-^{\beta,\nu}(t,x,\xi)\psi(t,x)$ has a compact support, say in $[0,T]\times[-R,R]\times[-R,R]$,
thus 
$\varphi_-(t,x,\xi)=g_-^{\beta,\nu}(t,x,\xi)\psi(t,x)$ is uniformly continuous
on this compact support.
Therefore for $\mu >0$, there exists $\gamma >0$ such that 
$|\varphi_-(t,x,\xi)-\varphi_-(0,x,\xi)| \leq \frac{\mu}{8R^2}$
for any $0\leq t<\gamma$ and any $x,\xi \in [-R,R]$,
and then for large $n$, we have $\eta_n < \gamma$ and
\begin{multline}
\left|\int_Q\int_\R f_+(t,x,\xi) \left(g_-^{\beta,\nu}(t,x,\xi)\psi(t,x)
-g_-^{\beta,\nu}(0,x,\xi)\psi(0,x)\right)
\omega_{\eta_n}'(t) d\xi dt dx \right|\\
\leq 
\int_Q\int_\R |f_+(t,x,\xi)| \rho_{\eta_n}(t) \frac{\mu}{8R^2} \1_{(x,\xi) \in [-R,R]^2}
d\xi dt dx\\
\leq \mu \int \rho_{\eta_n}(t) \,dt = \mu.
\label{eq:mutnul}\end{multline}
Thus we obtain, at the limit $n\to+\infty$ in (\ref{eqi}),
\begin{multline*}
\int_{Q}\int_\R -f_+g_-^{\beta,\nu}(\partial_t\psi+k(x) a(\xi)\partial_x\psi) d\xi dt dx\\
+\int_\R\int_\R -f_{0,+}g_-^{\beta,\nu}(0,x,\xi)\psi(0,x) d\xi dx \geq 0.
\end{multline*}

The next step is then to remove the hypothesis that $\psi$ vanishes at $x=0$ by setting $\psi(t,x)=\theta(t,x)\omega_{\eta_n}(x)$ where $\theta\in C^\infty_c([0,T[\times\R)$ is a non-negative test-function. We have
\begin{multline*}
\int_{Q}\int_\R -f_+g_-^{\beta,\nu}\omega_{\eta_n}(x)(\partial_t\theta
+k(x) a(\xi)\partial_x\theta) d\xi dt dx\\
+ \int_Q\int_\R -f_+g_-^{\beta,\nu}\theta(t,x) k(x) a(\xi)\omega_{\eta_n}'(x) d\xi dt dx\\
+\int_\R\int_\R -f_{0,+}g_-^{\beta,\nu}(0,x,\xi)\theta(0,x) \omega_{\eta_n}(x) d\xi dx \geq 0.
\end{multline*}

By \refe{tracespaceext} with $\theta_-(t,\xi)=g_-^{\beta,\nu}(t,0,\xi)\theta(t,0)$, 
\begin{multline*}
\lim_{n\to+\infty} \int_Q\int_\R f_+k(x) a(\xi)\omega_{\eta_n}'(x)
g_-^{\beta,\nu}(t,0,\xi)\theta(t,0)d\xi dt dx \\
\geq -(k_L-k_R)^+ \int_0^T\int_{\R}a(\xi) g_-^{\beta,\nu}(t,0,\xi)\theta(t,0) d\xi dt,
\end{multline*}
and by an argument similar to \refe{eq:mutnul}, the limit as $[n\to+\infty]$ of the term
\begin{equation*}
\int_Q\int_\R f_+k(x) a(\xi)\omega_{\eta_n}'(x)
 \left( 
 g_-^{\beta,\nu}(t,x,\xi)\theta(t,x)-
 g_-^{\beta,\nu}(t,0,\xi)\theta(t,0) \right)d\xi dt dx 
\end{equation*}
is zero. We have therefore
\begin{multline*}
\int_{Q}\int_\R -f_+g_-^{\beta,\nu}(\partial_t\theta+k(x) a(\xi)\partial_x\theta) d\xi dt dx\\
+(k_L-k_R)^+\int_0^T\int_\R a(\xi)  g_-^{\beta,\nu}(t,0,\xi)\theta(t,0) d\xi dt\\
+\int_\R\int_\R -f_{0,+}g_-^{\beta,\nu}(0,x,\xi)\theta(0,x) d\xi dx
 \geq 0.
\end{multline*}

Since $(k_L-k_R)^+=0$, we have actually
\begin{multline*}
\int_{Q}\int_\R -f_+g_-^{\beta,\nu}(\partial_t\theta+k(x) a(\xi)\partial_x\theta) d\xi dt dx\\
+\int_\R\int_\R -f_{0,+}g_-^{\beta,\nu}(0,x,\xi)\theta(0,x) d\xi dx\geq 0.
\end{multline*}
Take $\beta=\eta_n$ where $(\eta_n)$ is given in Prop.~\ref{prop:trace}. At the limit $\nu\to 0$ first, then $n\to+\infty$, we obtain
\begin{multline}
\int_{Q}\int_\R -f_+g_-(\partial_t\theta+k(x) a(\xi)\partial_x\theta) d\xi dt dx\\
+\limsup_{n\to+\infty}\int_\R\int_\R -f_{0,+}g_-^{\eta_n}(0,x,\xi)\theta(0,x) d\xi dx\geq 0.
\label{U4}\end{multline}
Observe that
\begin{align*}
g_-^{\eta_n}(0,x,\xi)=&\int_0^T g_-(t,x,\xi)\rho_{\eta_n}(t)dt\\
=&\int_0^T g_-(t,x,\xi)\omega_{\eta_n}'(t)dt.
\end{align*}
By \refe{tracetimeext} (transposed to $g_-$ tested against a function $\varphi_+$), we have
\begin{equation*}
\lim_{n\to+\infty}\int_\R\int_\R -f_{0,+}g_-^{\eta_n}(0,x,\xi)\theta(0,x) d\xi dx\leq \int_\R\int_\R -f_{0,+}g_{0,-}\theta(0,x) d\xi dx.
\end{equation*}
Since
\begin{equation*}
\int_\R -f_{0,+}g_{0,-} d\xi=\int_\R -\sgn_+(u_0-\xi)\sgn_-(v_0-\xi)d\xi=(u_0-v_0)^+,
\end{equation*}
we obtain by \refe{U4},
\begin{equation}
\int_Q\int_\R -f_+g_- (\partial_t\theta+k(x)a(\xi)\partial_x\theta) d\xi dt dx
 +\int_\R (u_0-v_0)^+\theta(0,x)dx\geq 0.
\label{U5}\end{equation}
It is then classical to conclude to \refe{eq:compfg}: let $M>0,R>MT$, let $\eta>0$ and let $r$ be a non-negative, non-increasing function such that $r\equiv 1$ on $[0,R]$, $r\equiv 0$ on $[R+\eta,+\infty[$. Set
$\theta(t,x)=\frac{T-t}{T}r(|x|+Mt)$ in \refe{U5} to obtain
\begin{equation*}
\frac{1}{T}\int_Q\int_\R -f_+g_- r(|x|+Mt)d\xi dt dx\leq \int_{\{|x|\leq R+\eta\}} (u_0-v_0)^+dx+\mathrm{J},
\end{equation*}
where the remainder term is
\begin{equation*}
\mathrm{J}=\int_Q\int_\R -f_+g_- \frac{T-t}{T}r'(|x|+Mt)
(M+k(x)a(\xi)\sgn(x)) d\xi dx dt.
\end{equation*}
By definition of $M$, $\mathrm{J}\leq 0$ and since $r(|x|+Mt)=1$ for $|x|\leq R-MT$, $0\leq t\leq T$, we obtain
\begin{equation*}
\frac{1}{T}\int_0^T\int_{|x|<R-MT}\int_\R -f_+g_- d\xi dx dt\leq \int_{\{|x|\leq R+\eta\}} (u_0-v_0)^+dx.
\end{equation*} 
Replacing $R$ by $R+MT$, and letting $\eta\to 0$ gives \refe{eq:compfg}. \qed

%%%%%%%%%%%%%%%%%%%%%%%%%%%%%%%%%%%%%%%%%%%%%%%%%%%%%%%%%%%%%%%%%%%%%%%%%%%%%%%%%%%
%%%%%%%%%%%%%%%%%%%%%%%%%%%%%%%%%%%%%%%%%%%%%%%%%%%%%%%%%%%%%%%%%%%%%%%%%%%%%%%%%%%
%%%%%%%%%%%%%%%%%%%%%%%%%%%%%%%%%%%%%%%%%%%%%%%%%%%%%%%%%%%%%%%%%%%%%%%%%%%%%%%%%%%
\section{Convergence of the BGK approximation}\label{sec:cvbgk}
%%%%%%%%%%%%%%%%%%%%%%%%%%%%%%%%%%%%%%%%%%%%%%%%%%%%%%%%%%%%%%%%%%%%%%%%%%%%%%%%%%%
%%%%%%%%%%%%%%%%%%%%%%%%%%%%%%%%%%%%%%%%%%%%%%%%%%%%%%%%%%%%%%%%%%%%%%%%%%%%%%%%%%%
%%%%%%%%%%%%%%%%%%%%%%%%%%%%%%%%%%%%%%%%%%%%%%%%%%%%%%%%%%%%%%%%%%%%%%%%%%%%%%%%%%%

\begin{theorem} Let $u_0\in L^1\cap L^\infty(\R)$, $0\leq u_0\leq 1$ a.e. When $\eps\to0$, the solution $f^\eps$ to the \refe{bgk} with initial datum $f_0=\chi_{u_0}$ converges in $L^p(Q\times\R_\xi)$, $1\leq p<+\infty$ to $\chi_u$, where $u$ is the unique solution to \refe{SCL}-\refe{IC}.
\label{th:cvbgk}\end{theorem}

{\bf Proof:} For $f\in L^1(\R_\xi)$, set 
\begin{equation*}
m_f(\xi)=\int_{-\infty}^\xi (\chi_u-f)(\zeta) d\zeta,\quad u=\int_\R f(\xi) d\xi.
\end{equation*}
It is easy to check that $m_f\geq 0$ if $0\leq \sgn(\xi)f(\xi)\leq 1$ for a.e. $\xi$ ({\it cf.} (29) in \cite{Brenier83}). In our context, we have $0\leq f^\eps\leq \chi_1$, hence $m^\eps:=\frac{1}{\eps}m_{f^\eps}\geq 0$. Viewed as a measure, $m^\eps$ is supported in $[0,T]\times\R_x\times[0,1]$. Integration with respect to $\xi$ in \refe{bgk} gives
\begin{equation*}
m^\eps(\xi)=\partial_t\left(\int_0^\xi f^\eps(\zeta)d\zeta\right)+\partial_x\left( k(x)\int_0^\xi a(\zeta)f^\eps(\zeta)d\zeta\right)
\end{equation*}
in $\mathcal{D}'(]0,T[\times\R_x)$. Summing over $(t,x)\in [0,T]\times [x_1,x_2]$, $\xi\in ]0,1[$, we get the estimate
\begin{multline}
m^\eps([0,T]\times[x_1,x_2]\times[0,1])=\int_{x_1}^{x_2}\int_0^1(1-\xi) 
(f^\eps(T,x,\xi)-f^\eps(0,x,\xi))d\xi dx\\
	+\left[\int_0^T\int_0^1 (1-\xi)k(x)a(\xi)f^\eps(t,x,\xi) d\xi dt \right]_{x_1}^{x_2}.
\label{cv1}\end{multline} 
Since $f^\eps(t)\in L^1(\R_x\times\R_\xi)$, there exists sequences $(x_1^n)\downarrow -\infty$ and $(x_2^n)\uparrow+\infty$ such that the last term of the right hand-side in \refe{cv1} tends to $0$ when $n\to+\infty$. Since, besides, $f^\eps\geq 0$ and 
\begin{equation*}
\int_\R\int_\R f^\eps(T,x,\xi) d\xi dx \leq \int_\R\int_\R \chi_{u_0} d\xi dx=\|u_0\|_{L^1(\R)},
\end{equation*}
we obtain the uniform estimate
\begin{equation}
m^\eps([0,T]\times\R\times[0,1])\leq\|u_0\|_{L^1(\R)}.
\label{cv2}\end{equation}
We also have
\begin{equation}
0\leq f^\eps\leq\chi_1,\quad -\partial_\xi f^\eps_+(t,x,\xi)=\nu^\eps_{t,x}(\xi)
+\mathcal{O}(\eps)
\label{cv3}\end{equation}
where $\nu^\eps_{t,x}(\xi):=\delta_{u^\eps(t,x)}(\xi)$ and the identity is satisfied in $\mathcal{D}'(]0,T[\times\R_x\times\R_\xi)$. 
Indeed, by \refe{bgk},
\begin{equation*}
f^\eps=\chi_{u^\eps}+\eps(\partial_t f^\eps+\partial_x(k(x)a(\xi)f^\eps))=\chi_{u^\eps}+\mathcal{O}(\eps),
\end{equation*}
hence 
\begin{equation*}
-\partial_\xi f^\eps_+=-\partial_\xi f^\eps+\delta_0(\xi)=-\partial_\xi\chi_{u^\eps}+\delta_0(\xi)+\mathcal{O}(\eps)=\delta_{u^\eps}(\xi)
+\mathcal{O}(\eps).
\end{equation*}
Notice that, for a.e. $(t,x)$, $\nu^\eps_{t,x}$ is supported in the fixed compact subset $[0,1]$ of $\R_\xi$. We deduce from \refe{cv2}-\refe{cv3} that, up to a subsequence, there exists a non-negative measure $m$ on $\R^3$ supported in $[0,T]\times\R_x\times[0,1]$, a function $f\in L^\infty(]0,T[;L^1(\R_x\times\R_\xi))$ such that $0\leq f\leq\chi_1$, $-\partial_\xi f_+(t,x,\xi)=\nu_{t,x}(\xi)$ where $\nu$ is a Young measure $Q\to\R_\xi$ and such that $m^\eps\rightharpoonup m$ weakly in the sense of measures ({\it i.e.} $\<m^\eps-m,\varphi\>\to 0$ for every continuous compactly supported $\varphi$ on $\R^3$) and $f^\eps\rightharpoonup f$ in $L^\infty(Q\times\R_\xi)$ weak-star. Besides, since $f^\eps\in C([0,T];L^1_{x,\xi})$ satisfies $f^\eps(0)=f_0$ and the BGK equation
\begin{equation*}
\partial_t f^\eps+\partial_x(k(x)a(\xi)f^\eps)=\partial_\xi m^\eps,
\end{equation*}
it satisfies the weak formulation: for all $\psi\in C^\infty_c([0,T[\times\R\times\R)$,
\begin{multline*}
\int_{Q}\int_\R f^\eps (\partial_t\psi+k(x) a(\xi)\partial_x\psi) d\xi dt dx
+\int_{\R}\int_\R f_{0}\psi(0,x,\xi) d\xi dx\\
=\int_{Q}\int_\R\partial_\xi\psi dm^\eps(t,x,\xi).
\end{multline*}
In particular, we have
\begin{align}
&\int_{Q}\int_\R f^\eps_\pm (\partial_t\psi+k(x) a(\xi)\partial_x\psi) d\xi dt dx
+\int_{\R}\int_\R f_{0,\pm}\psi(0,x,\xi) d\xi dx\nonumber\\
=&-\int_{Q}\int_\R\sgn_\mp(\xi)k(x)a(\xi)\partial_x\psi d\xi dt dx+\int_{Q}\int_\R\partial_\xi\psi dm^\eps(t,x,\xi)\nonumber\\
=&(k_R-k_L)\int_0^T\int_\R\sgn_\mp(\xi)a(\xi)\psi(t,0,\xi) d\xi dt+\int_{Q}\int_\R\partial_\xi\psi dm^\eps(t,x,\xi)\nonumber\\
=&(k_L-k_R)^\pm\int_0^T\int_\R a(\xi)\psi(t,0,\xi)d\xi dt+\int_{Q}\int_\R\partial_\xi\psi dm^\eps_\pm(t,x,\xi),\label{weakgeseps}
\end{align}
where 
\begin{multline}
\<m^\eps_\pm,\partial_\xi\psi\>:=\<m^\eps,\partial_\xi\psi\>\\
-\int_0^T\int_\R a(\xi)[(k_L-k_R)\sgn_\mp(\xi)+(k_L-k_R)^\pm]\psi(t,0,\xi) d\xi dt.
\label{cv5}\end{multline}
More precisely, we set 
\begin{equation*}
m_+^\eps=m^\eps+\int_\xi^{+\infty}a(\zeta)[(k_L-k_R)^+\sgn_+(\zeta)-(k_L-k_R)^-\sgn_-(\zeta)]d\zeta
\delta(x=0),
\end{equation*}
and
\begin{equation*}
m_-^\eps=m^\eps+\int_{-\infty}^\xi a(\zeta)[(k_L-k_R)^+\sgn_+(\zeta)-(k_L-k_R)^-\sgn_-(\zeta)]d\zeta
\delta(x=0).
\end{equation*}
Notice that in both cases, and since $A(\xi)\geq 0$ for any $\xi$, we have added a non-negative quantity to $m^\eps$. At the limit $\eps\to 0$ we thus obtain $m^\eps_\pm\rightharpoonup m_\pm$ where $m_\pm$ is a non-negative measure. Examination of the support of $m^\eps_\pm$ shows that $m_+$, {\it resp. $m_-$} is supported in $[0,T]\times\R_x\times]-\infty,1]$, {\it resp.} $[0,T]\times\R_x\times[0,+\infty[$. At the limit $\eps\to 0$, we thus obtain the kinetic formulation \refe{eq:kineticf}. We conclude that $f$ is a generalized solution to \refe{SCL}-\refe{IC}. By Theorem~\ref{th:redU}, $f=\chi_u$ where $u\in L^\infty(Q)$ is solution to \refe{SCL}-\refe{IC}. By uniqueness, the whole sequence $(f^\eps)$ converges (in $L^\infty$ weak-star) to $\chi_u$. Actually the convergence is strong since
\begin{align}
\int_Q\int_\R |f^\eps-\chi_u|^2 d\xi dt dx=&\int_Q\int_\R |f^\eps|^2-2f^\eps\chi_u+\chi_u d\xi dt dx\nonumber\\
\leq &\int_Q\int_\R f^\eps-2f^\eps\chi_u+\chi_u d\xi dt dx.\label{cv6}
\end{align}
We have used the fact that $0\leq f^\eps\leq 1$. The right-hand side of \refe{cv6} tends to $0$ when $\eps\to 0$ since $1$, $\chi_u \in L^\infty$ can be taken as test functions.
Hence $f^\eps\to \chi_u$ in $L^2(Q\times\R)$. The convergence in $L^p(Q\times\R)$, $1\leq p<+\infty$ follows from the uniform bound on $f^\eps$ in $L^1\cap L^\infty(Q\times\R)$. \qed
\bigskip

{\bf Remark:} it is possible to relax the assumption that the initial datum for \refe{bgk} is at equilibrium and independent on $\eps$ in Theorem~\ref{th:cvbgk}. Indeed, the conclusion of Theorem~\ref{th:cvbgk} remains valid under the hypothesis that the initial datum $f_0^\eps$ for \refe{bgk} satisfies
\begin{equation}
0\leq f^\eps_0\leq\chi_1,\quad f^\eps_0\rightharpoonup f_0,\quad u_0(x):=\int_\R f_0(x,\xi)d\xi,
\label{iceps}\end{equation}
where $f^\eps_0\rightharpoonup f_0$ in \refe{iceps} denotes weak convergence in $L^1(\R_x\times\R_\xi)$. Indeed, the proof of Theorem~\ref{th:cvbgk} remains unchanged under the following modification: passing to the limit in \refe{weakgeseps}, we obtain that $f$ is a generalized solution to \refe{SCL} with an initial datum $f_0$ that is not necessary at equilibrium. However, we have ({\it cf.} (29) in \cite{Brenier83})
\begin{equation*}
f_0-\sgn_\mp(\xi)=\sgn_\pm(u_0-\xi)-\partial_\xi m^0_\pm,
\end{equation*}
where $m^0_+$ ({\it resp.} $m^0_-$) is a non-negative measure supported in $[0,T]\times\R\times]-\infty,1]$ ({\it resp.} $[0,T]\times\R\times[0,+\infty[$). Consequently, up to a modification of the kinetic measure $m_\pm$, we obtain that $f$ is indeed a generalized solution to \refe{SCL}-\refe{IC}. The rest of the proof is similar.

\providecommand{\bysame}{\leavevmode\hbox to3em{\hrulefill}\thinspace}
\providecommand{\MR}{\relax\ifhmode\unskip\space\fi MR }
\providecommand{\MRhref}[2]{%
  \href{http://www.ams.org/mathscinet-getitem?mr=#1}{#2}
}
\providecommand{\href}[2]{#2}

\end{document}